\documentclass[11pt, oneside]{article}   	
\usepackage{geometry}                		
\usepackage{graphicx}				
\usepackage{amssymb}

\usepackage{theorem}
\usepackage[colorlinks=true]{hyperref}
\theoremheaderfont{\normalfont\bfseries}
\theorembodyfont{\itshape}

\newtheorem{theorem}[subsection]{Theorem}
\newtheorem{proposition}[subsection]{Proposition}
\newtheorem{lemma}[subsection]{Lemma}

\theorembodyfont{\upshape}
{\theoremstyle{break}

}
\theorembodyfont{\upshape}
{\theoremstyle{break}
\newtheorem{remark}[subsection]{Remark}
}
\newcommand{\qed}{\nopagebreak\par\hspace*{\fill}$\square$\par\vskip2mm}

\newcommand{\Um}{\mathop{\mathrm{Um}}\nolimits}
\newcommand{\SL}{\mathop{\mathit{SL}}\nolimits}
\newcommand{\WMS}{\mathop{\mathit{WMS}}\nolimits}
\newcommand{\V}{{\Um_{\mathrm{ gen}}(R)}}
\newcommand{\ba}{{\mathbf a}}
\newcommand{\bb}{{\mathbf b}}
\title{Extrapolating an Euler class}
\author{Wilberd van der Kallen}
\date{}							

\begin{document}
\maketitle
\sloppy
\begin{abstract}
Let $R$ be a noetherian ring of dimension $d$ and let $n$ be an integer so that $n\leq d\leq 2n-3$.
Let $(a_1,\dots,a_{n+1})$ be a unimodular row so that the ideal $J=(a_1,\dots,a_n)$ has height $n$.
Jean Fasel has associated to this row an element $[(J,\omega_J)]$ in the Euler class group $E^n(R)$,
with $\omega_J:(R/J)^n\to J/J^2$ given by $(\bar a_1,\dots,\bar a_{n-1},\bar a_n\bar a_{n+1})$.
If $R$ contains an infinite field $F$ then we show that the rule of Fasel defines a homomorphism
from $\WMS_{n+1}(R)=\Um_{n+1}(R)/E_{n+1}(R)$ to $E^n(R)$. The main problem is to get a well
defined map on all of $\Um_{n+1}(R)$. Similar results have been obtained by 
Das and Zinna \cite{DasZinna}, with a different proof. Our proof uses that every Zariski  open subset of $\SL_{n+1}(F)$
is path connected for walks made up of elementary matrices. 
\end{abstract}
\section{Recollections}
\subsection{The group of orbits}Let $n\geq3$.
Let $R$ be a commutative noetherian ring of Krull dimension $d$, $d\leq2n-2$. 
As usual $E_m(R)$ denotes the subgroup of $\SL_m(R)$ generated by elementary matrices $e_{ij}(r)$
and $\Um_m(R)$ denotes the set of unimodular rows of length $m$ over $R$.
Then 
\cite[Theorem 4.1]{vdk module} provides an abelian group structure on the orbit set
$\Um_{n+1}(R)/E_{n+1}(R)$.
The abelian group that is obtained is called $\WMS_{n+1}(R)$.
As explained in \cite[\S3]{vdk mennicke} the group law may be characterized as follows.
 If $\alpha, \beta\in \WMS_{n+1}(R)$, one may choose representatives $(a_1,\dots,a_{n+1})\in\alpha$,
$(b_1,\dots,b_{n+1})\in\beta$ so that $a_1+b_1=1$ and $a_i=b_i$ for $i>1$.
Then $(a_1b_1,a_2,\dots,a_{n+1})$ is a representative of $\alpha+\beta$.
This rule reflects the homotopic join of Borsuk \cite{Borsuk}.
\subsection{The Euler class group}
From now on let $3\leq d\leq2n-3$.
In \cite{BhatwadekarSridharan2} the authors introduce an Euler class group $E^n(R)$ generalizing the Euler class group of 
\cite{BhatwadekarSridharan1}. The latter  corresponds with the case $n=d$.
The Euler class group is an abelian group given by a presentation. Generators are pairs $(J,\omega_J)$ where $J$ is a height $n$ ideal in $R$ equipped with
a surjective map $(R/J)^n\to J/J^2$. Think of a codimension $n$ subvariety with trivial conormal bundle together with a trivialization of said bundle.
Relations are
\subsubsection*{Disconnected sum}Let $(J,\omega_J)$ be a generator.
If $J=KL$ with $K$, $L$ comaximal ideals of height $n$, then $R/J=R/K\times R/L$ and
 $J/J^2=K/K^2\times L/L^2$, so that $\omega_J=\omega_k\times \omega_L$.
The relation is $$(J,\omega_J)=(K,\omega_K)+(L,\omega_L).$$ 
\subsubsection*{Complete intersection}
Let $(J,\omega_J)$ be a generator such that $\omega_J$ lifts to a surjection $R^n\to J$.
Then $$(J,\omega_J)=0.$$ 
\subsubsection*{Elementary action}
Let $(J,\omega_J)$ be a generator and let $g\in E_n(R/J)$.
Then $$(J,\omega_J)=(J,\omega_J\circ  g).$$
One may define $E^n(R)$ by taking the disconnected sum relations and the complete intersection relations as defining relations, cf.\
\cite[Proposition 2.2]{DasZinna}. We denote 
the class of $(J,\omega_J)$ in $E^n(R)$ by $[(J,\omega_J)]$. One shows with \cite[Corollary 2.4, Proposition 3.1]{BhatwadekarSridharan2}
that every element of $E^n(R)$ can be written in the form
$[(J,\omega_J)]$. And one shows as in \cite[Proposition 2.2]{DasZinna} that the elementary action relations also hold. (Use \cite[Corollary 2.4]{BhatwadekarSridharan2} and use that
$E_n(R)\to E_n(R/K)\times E_n(R/L)$ is surjective in the disconnected sum setting.) Note that  \cite[Corollary 2.4]{BhatwadekarSridharan2} 
only needs $d\leq2n-1$.

\subsection{The old homomorphism and the new one}
If our $R$ is a regular ring containing an infinite field then Bhatwadekar and Sridharan   define a homomorphism $\WMS_{n+1}(R)\to E^n(R)$ with useful properties when $n$ is even (\cite[Theorem 5.7]{BhatwadekarSridharan2}).
But it vanishes when $n$ is odd. Jean Fasel noticed that in $\mathbb A^1$ homotopy one can do better. He proposed a formula that would
also be useful when $n$ is odd. In fact the same formula was discussed by Bhatwadekar and Sridharan after \cite[Theorem 7.3]{BhatwadekarSridharan1},
in the case of even $n$.
It is already known that 
the formula of Fasel works for $3\leq n=d$ \cite[Theorem 3.6]{DasZinna}. That is, it defines a homomorphism $\WMS_{n+1}(R)\to E^n(R)$.
If $R$ is a domain it is also known to work
\cite[Remark 3.11]{DasZinna} (always assuming $3\leq d\leq2n-3$). 
Our purpose is to show that his formula works when $R$ contains an infinite field $F$. The main difference between this note and \cite{DasZinna} is in
the proof strategy. Rather than studying $E^n(R)$ more closely, as is done
in \cite{DasZinna}, we concentrate on $\Um_{n+1}(R)$.
We use paths made up of elementary matrices in $\SL_{n+1}(F)$ to walk back and forth between general unimodular rows and rows for which 
we already know what to do.

\section{Elementary paths}
The group $E_{n+1}(R)$ is generated by elementary matrices $e_{ij}(r)$.
We call a sequence $g_1$, \dots, $g_m$ of elements of $E_{n+1}(R)$ an \emph{elementary path} if the $g_i^{-1}g_{i+1}$ are elementary for
$i=1,\dots,m-1$. We call it an $F$-path if moreover all $g_i$ are in $\SL_{n+1}(F)$.
Notice that if $g_1$, \dots, $g_m$ is a path, then the reverse sequence $g_m$, \dots, $g_1$ is also a path.

We provide $\SL_{n+1}(F)$ with the topology induced by the Zariski topology on the algebraic group $\SL_{n+1}$ defined over $F$.
We say that a subset $U$ of $\SL_{n+1}(F)$ is path connected if any two elements of $U$ can be joined by an $F$-path that stays within $U$.

\begin{proposition}\label{omega}
Any nonempty open subset of $\SL_{n+1}(F)$ is path connected.
\end{proposition}
\paragraph{Proof}We give two proofs.
Let $\Omega$ be the big cell of $\SL_{n+1}(F)$. Then by \cite[Prop 2.6]{vdk units} every open subset of $\Omega$ is path connected. But then for any 
$g\in \SL_{n+1}(F)$ every open subset of $g\Omega$ is path connected. So if $U$ is a nonempty open subset of $\SL_{n+1}(F)$, then it is covered by
mutually intersecting path connected subsets of type $U\cap g\Omega$.

For the second proof recall that an element $g\in \SL_{n+1}(F)$ that is in general position may be reduced to the identity matrix 
with just $N=(n+1)^2-1$ elementary operations. It goes like this. Add a multiple of the last column to the first to achieve that $g_{11}$ becomes equal to
one. Then clear
the first row in $n$ steps. Add a multiple of the last column to the second to achieve that $g_{22}$ becomes equal to one. Then clear
the second row in $n$ steps. Keep going until the second last row has been cleared. Notice that $g_{n+1,n+1}$ has become equal to one. 
Clear the last row
in $n$ steps. Reversing this procedure one finds a sequence $\alpha_1$, \dots, $\alpha_{N}$ of roots so that the map
$(t_1,\dots,t_{N})\mapsto e_{\alpha_1}(t_1)\cdots e_{\alpha_N}(t_N)$ defines a birational map from $F^N$ to $\SL_{n+1}(F)$.
Now if $p$, $q$ are elements of the open subset $U$, then the condition that the path $p$, $pe_{\alpha_1}(t_1)$, 
$pe_{\alpha_1}(t_1)e_{\alpha_2}(t_2)$, \dots\ stays inside $U$ defines an open condition on $F^N$. So there is a nonempty open subset
of $U$ of elements that can be reached by an $F$-path within $U$ that starts at $p$. Similarly there is a nonempty open subset
of $U$ that can be reached by an $F$-path starting at $q$. These two subsets intersect.
\qed

\begin{remark}
Of course there is no such result for finite fields.
\end{remark}

\subsection{Prime avoidance}
We will tacitly use variations on the proof of the following classical lemma.

\begin{lemma}
Let $(a_1,\dots,a_{m+1})$ be a unimodular row over $R$.
There are $\lambda_i\in R$ so that the ideal $(a_1+\lambda_1 a_{m+1},\dots,a_m+\lambda_ma_{m+1})$ has height $m$ or equals $R$.
\end{lemma}
\paragraph{Proof}
We argue by induction on $m$. For $m=1$ it is a prime avoidance exercise \cite[Exercise 3.19]{Eisenbud}. Let $m>1$.
The set of ideals of the form $(a_1+\lambda_1 a_{m+1},\dots,a_m+\lambda_ma_{m+1})$ does not change when we add multiples of
$a_2$, \dots, $a_{m+1}$ to $a_1$. So by the same exercise we may assume $a_1$ avoids all minimal primes. Apply the inductive hypothesis to the 
unimodular row $(\bar a_2,\dots,\bar a_{m+1})$ over $R/(a_1)$.
\qed

\subsection{Generic unimodular rows}

We define $\V$ to be the set of unimodular rows $\ba=(a_1,\dots,a_{n+1})$ over $R$ for which the ideal $(a_na_{n+1})$ has height one.
We call such rows generic. So in  a generic row the last two entries avoid every minimal prime ideal of $R$.

If $\ba\in \Um_{n+1}(R)$ and $g_1$, \dots, $g_m$ is an elementary path in $E_{n+1}(R)$, then we call the sequence
$\ba g_1$, \dots, $\ba g_m$ a 
path in $ \Um_{n+1}(R)$.
If moreover $g_1$, \dots, $g_m$ is an $F$-path, then we also call $\ba g_1$, \dots, $\ba g_m$ an $F$-path.

\begin{proposition}Generic rows detect orbits:

\begin{enumerate}
\item Every $\SL_{n+1}(F)$-orbit in $\Um_{n+1}(R)$ intersects $\V$ in a nonempty path connected subset.\label{F intersect}
\item Every $E_{n+1}(R)$-orbit in $\Um_{n+1}(R)$ intersects $\V$ in a nonempty path connected subset.\label{R intersect}
\end{enumerate}
\end{proposition}
\paragraph{Proof}
For $\ba\in \Um_{n+1}(R)$ consider the set of $g\in \SL_{n+1}(F)$ for which $\ba g\in\V$. It is a nonempty open subset of $\SL_{n+1}(F)$.
(The complement is closed because minimal primes are linear subspaces.) Therefore part \ref{F intersect} follows from Proposition \ref{omega}.

To prove the second part, fix a path component $P$ of $\V$ and let $X$ be the set of $\ba\in  \Um_{n+1}(R)$ for which there is an $F$-path
starting at $\ba$ and ending in $P$. Clearly $X$ is invariant under the action by $\SL_{n+1}(F)$. Notice that by part \ref{F intersect},
if $\ba\in X$, then every $F$-path from 
$\ba$ to $\V$ lands in $P$.
If we show that $X$ is also invariant under the action by $e_{21}(r)$ for $r\in R$, then it will follow that $X$ is an $E_{n+1}(R)$-orbit
with $X\cap\V=P$. So let $\ba\in X$. We need to show that $\ba e_{21}(r)\in X$. We may replace $\ba$ with $\ba g$ for any $g\in \SL_{n+1}(F)$ 
that commutes with $e_{21}(r)$. Therefore we may assume $\ba=(a_1,\dots,a_{n+1})$ is such that there are $\lambda$, $\mu\in F$ with
$\ba e_{1,n}(\lambda)e_{1,n+1}(\mu)\in\V$ and $\ba e_{21}(r) e_{1,n}(\lambda)e_{1,n+1}(\mu)\in\V$. Note that $\ba e_{1,n}(\lambda)e_{1,n+1}(\mu)\in P$.
We can walk inside $\V$ from $\ba e_{21}(r) e_{1n}(\lambda)e_{n+1}(\mu)$ to $\ba e_{1n}(\lambda)e_{n+1}(\mu)$ by way of
$(a_1,a_2,\dots,a_{n-1},a_n+\lambda (a_1+ra_2),a_{n+1}+\mu (a_1+ra_2))$ and
$(a_1,a_2,\dots,a_{n-1},a_n+\lambda a_1,a_{n+1}+\mu (a_1+ra_2))$.
\qed

\section{The map}
Before defining a map $\phi:\Um_{n+1}(R)\to E^n(R)$ we will define one on generic rows.
But first we define $\phi_0(\ba)$ when $\ba=(a_1,\dots,a_{n+1})$ is such  that the ideal $J=(a_,\dots,a_n)$ has height at least $n$.
If $J=R$ we put $\phi_0(\ba)=0$. Remains the case that $J$ has height $n$.
Then we follow Fasel and put 
$$\phi_0(\ba)=[(J,\omega_J)]$$
with $\omega_J:(R/J)^n\to J/J^2$ given by $(\bar a_1,\dots,\bar a_{n-1},\bar a_n\bar a_{n+1})$.

Now let $\ba=(a_1,\dots,a_{n+1})\in\V$ be generic.
Choose $\lambda_i$, $\mu_i$ in $R$ so  that $J_1=(a_1+\lambda_1a_n,\dots,a_{n-1}+\lambda_{n-1}a_n,a_{n+1})$ and
 $J_2=(a_1+\mu_1a_{n+1},\dots,a_{n-1}+\mu_{n-1}a_{n+1},a_{n})$ have height at least $n$.

\begin{lemma}
The sum of  $$[(J_1,\omega_{J_1})]:=\phi_0((a_1+\lambda_1a_n,\dots,a_{n-1}+\lambda_{n-1}a_n,a_{n+1},a_n))$$ and
 $$[(J_2,\omega_{J_2})]:=\phi_0((a_1+\mu_1a_{n+1},\dots,a_{n-1}+\mu_{n-1}a_{n+1},a_{n},a_{n+1}))$$ vanishes.

Therefore $[(J_1,\omega_{J_1})]$ does not depend on the choice of the $\lambda_i$ and 
$[(J_2,\omega_{J_2})]$ does not depend on the choice of the $\mu_i$.
\end{lemma}
\paragraph{Proof}(Compare \cite[Proposition 3.4]{DasZinna}.) As $\bar a_n$ is invertible in $R/J_1$ we get 
that
$[(J_1,\omega_{J_1})]$ equals 
$$\phi_0((a_1+\lambda_1a_n+\mu_1a_{n+1},\dots,a_{n-1}+\lambda_{n-1}a_n+\mu_{n-1}a_{n+1},a_{n+1},a_{n})).$$
Similarly $[(J_2,\omega_{J_2})]$ equals
$$\phi_0((a_1+\lambda_1a_n+\mu_1a_{n+1},\dots,a_{n-1}+\lambda_{n-1}a_n+\mu_{n-1}a_{n+1},a_{n},a_{n+1})).$$
As $J_1$, $J_2$ are comaximal, $[(J_1,\omega_{J_1})]+[(J_2,\omega_{J_2})]$  equals
$$
\phi_0((a_1+\lambda_1a_n+\mu_1a_{n+1},\dots,a_{n-1}+\lambda_{n-1}a_n+\mu_{n-1}a_{n+1},a_{n+1}a_n,1)).$$
\qed

\subsection{The map on generic rows}
For $\ba=(a_1,\dots,a_{n+1})\in\V$ we choose the $\mu_i$ as above and put 
$$\phi(\ba)=\phi_0((a_1+\mu_1a_{n+1},\dots,a_{n-1}+\mu_{n-1}a_{n+1},a_{n},a_{n+1})).$$
Note that $\phi((a_1,\dots,a_{n+1}))= -\phi((a_1,\dots,a_{n-1},a_{n+1},a_n))$.
Note also that, if $i\neq j$, $j<n$, then $\phi(\ba)$ does not change if we add a multiple of $a_i$ to $a_j$.

\begin{proposition}
The map $\phi$ is constant on path components of $\V$.
\end{proposition}
\paragraph{Proof}We must show that $\phi(\ba)=\phi(\ba e_{ij}(r))$
if $\ba=(a_1,\dots,a_{n+1})\in\V$ and  $e_{ij}(r)$ is an elementary matrix so that $ \ba e_{ij}(r)\in\V$.
We already know it when $j<n$. Remain cases with $j=n$ or $j=n+1$.

For instance, let $i=n$, $j=n+1$. We may add multiples of $a_n$, $a_{n+1}$ to $a_1$ through $a_{n-1}$ to reduce to the case that 
$\phi_0(\ba)$ is defined. But then $\phi_0(\ba)$ and $\phi_0(\ba e_{ij}(r))$ are computed as the class of the same $(J,\omega_J)$.

Next try $i=n-1$, $j=n+1$. Choose $\lambda\in R$ so that $(a_{n-1}+\lambda a_n)$ and $(a_{n+1}+r(a_{n-1}+\lambda a_n))$ both have height one.
Then $\phi(\ba)=\phi((a_1,\dots,a_{n-2},a_{n-1}+\lambda a_n,a_n,a_{n+1}))$ and 
$\phi(\ba e_{n-1,n+1}(r))=\phi((a_1,\dots,a_{n-2},a_{n-1}+\lambda a_n,a_n,a_{n+1}+r(a_{n-1}+\lambda a_n)))$ by the cases that have already
 been treated. So replacing $a_{n-1}$ with $a_{n-1}+\lambda a_n$ we may further assume that $(a_{n-1})$ has height one.
Adding multiples of $a_{n-1}$, $a_n$, $a_{n+1}$ to $a_1$ through $a_{n-2}$ we may further arrange that $(a_1,\dots, a_{n-1})$ has height $n-1$.
And adding a multiple of $a_n$ to $a_{n+1}$ we may further assume that $(a_1,\dots, a_{n-1},a_{n+1})$ and
$(a_1,\dots, a_{n-1},a_{n+1}+ra_{n-1})$ have height $n$.
Choose $\mu$ so that $(a_n+\mu a_{n+1})$ has height one and so that 
$\phi(\ba e_{n-1,n+1}(r))=\phi_0((a_1,\dots,a_{n-1},a_n+\mu(a_{n+1}+ra_{n-1}),a_{n+1}+ra_{n-1}))$.
Then $\phi(\ba e_{n-1,n+1}(r))=\phi_0((a_1,\dots,a_{n-1},a_n+\mu(a_{n+1}+ra_{n-1}),a_{n+1}))$, which equals
$-\phi_0((a_1,\dots,a_{n-1},a_{n+1},a_n+\mu(a_{n+1}+ra_{n-1})))$.
Therefore $\phi(\ba e_{n-1,n+1}(r))=-\phi_0((a_1,\dots,a_{n-1},a_{n+1},a_n+\mu a_{n+1}))=\phi(\ba)$.

Other cases are analogous or easier.
\qed

It is now clear how to extend $\phi$ to all of $\Um_{n+1}(R)$.
For $\ba\in\Um_{n+1}(R)$ take a path to $\V$ and define $\phi(\ba)$ to be the value of $\phi$ at the end of the path.
Note that $\phi$ now also extends the map $\phi_0$.

\begin{theorem}
Let $3\leq d\leq2n-3$ and let $R$ be a commutative noetherian ring of Krull dimension $d$.
If $R$ contains an infinite field then $\phi_0$ induces a homomorphism $\WMS_{n+1}(R)\to E^n(R)$.
\end{theorem}
\paragraph{Proof}To see that $\phi$ defines a homomorphism $\WMS_{n+1}(R)\to E^n(R)$ consider 
$\ba=(a_1,\dots,a_{n+1})$,
$\bb=(b_1,\dots,b_{n+1})$ so that $a_1+b_1=1$ and $a_i=b_i$ for $i>1$.
We may add multiples of $a_{n+1}$ to the other entries to arrange that moreover
$\phi_0(\ba)$, $\phi_0(\bb)$ are defined. 
Now use that the ideals $(a_1,\dots,a_{n})$,
$(b_1,\dots,b_{n})$ are comaximal.
\qed


\begin{thebibliography}{99}
\bibitem[B]{Borsuk}K. Borsuk, \emph{Remarks on the homotopic join of
maps},  Fundamenta Mathematicae {\bf (50)} 1961/62, 195--206.
\bibitem[BS1]{BhatwadekarSridharan1}S.M. Bhatwadekar and Raja Sridharan,
\emph{Euler
class group of a noetherian ring}, Compositio Math.\ {\bf (122)} 2000,
183--222.
\bibitem[BS2]{BhatwadekarSridharan2}S.M. Bhatwadekar and Raja Sridharan,
\emph{On Euler Classes and Stably Free Projective Modules}, 
Proceedings of the International Colloquium on Algebra, Arithmetic and Geometry 
Mumbai 2000, Part I, 139--158. Narosa Publishing House. International distribution by American Math Society.
\bibitem[DZ]{DasZinna}Mrinal Kanti Das and Md. Ali Zinna, \emph{``Strong" Euler class of a stably free module of odd rank}, Journal of Algebra {\bf 432} (2015), 185--204.
\bibitem[Eis]{Eisenbud}D.~Eisenbud,
\emph{Commutative algebra.
With a view toward algebraic geometry}. Graduate Texts in Mathematics, 150.
Springer-Verlag, New York, 1995.
\bibitem[vdK1]{vdk units}W. van der Kallen,
\emph{The $K_2$ of rings with many units},
Ann. Sci. \'Ec. Norm. Sup. (4) 10 (1977), 473--515.
\bibitem[vdK2]{vdk module}W. van der Kallen,  \emph{A module structure on
certain
orbit
sets of unimodular rows},
 Journal  of Pure and Appl.\ Algebra {\bf (57)} 1989, 281--316.
 \bibitem[vdK3]{vdk mennicke}
 W. van der Kallen, \emph{From Mennicke symbols to Euler class groups}, 
Proceedings of the International Colloquium on Algebra, Arithmetic and Geometry 
Mumbai 2000, Part II, 341--354. Narosa Publishing House. International distribution by American Math Society.
\end{thebibliography}
\end{document}